\documentclass[12pt,reqno]{amsart}

\usepackage{amssymb,amsmath, amsfonts, latexsym}
\usepackage{enumerate}
\newtheorem{thm}{Theorem}[section]
\newtheorem{cor}[thm]{Corollary}
\newtheorem{lem}[thm]{Lemma}
\newtheorem{prop}[thm]{Proposition}

\newtheorem{rem}{Remark}[section]

\newcommand{\dR}{\mathbb{R}}
\newcommand{\dP}{\mathbb{P}}

\newcommand{\superexp} {\quad\underset{\lambda_n}{\overset{\rm superexp}{\longrightarrow}} }

\newcommand{\superexpldp} {\quad\underset{n}{\overset{\rm superexp}{\longrightarrow}} }
\newcommand{\superexpmdp} {\quad \underset{v^2_n}{\overset{\rm superexp}{\longrightarrow}} }

\def\build#1_#2^#3{\mathrel{\mathop{\kern 0pt#1}\limits_{#2}^{#3}}}

 \numberwithin{equation}{section}

\begin{document}

\title[LDP of threshold estimator]{Large deviations  of the Threshold estimator of integrated (co-)volatility vector in the presence of jumps}

\author{Hac\`ene Djellout}
\email{Hacene.Djellout@math.univ-bpclermont.fr}
\address{Laboratoire de Math\'ematiques, CNRS UMR 6620, Universit\'e Blaise Pascal, Avenue des Landais,BP80026, 63171 Aubi\`ere Cedex, France.}

\author{Hui  Jiang}
\email{huijiang@nuaa.edu.cn}
\address{Department of Mathematics, Nanjing University of Aeronautics and Astronautics,
29 Yudao Street, Nanjing 210016, China.}


\date{\today}
\begin{abstract} Recently a considerable interest has been paid on the estimation problem of the realized volatility and covolatility  by using high-frequency data of financial price processes in financial econometrics.  Threshold estimation is one of the useful techniques in the inference for jump-type stochastic processes from discrete observations. In this paper, we adopt the threshold estimator introduced by Mancini \cite{Mancini4} where only the variations under a given threshold function are taken into account. The purpose of this work is to investigate large and moderate deviations for the threshold estimator of the integrated variance-covariance vector. This paper is an extension of the previous work in Djellout et al \cite{Djellout3}.  where the problem has been studied in absence of the jump component. We will use the approximation lemma to prove the LDP. As the reader can expect we obtain the same results as in the case without jump.
\end{abstract}

\maketitle

\vspace{-0.5cm}

\begin{center}
\end{center}

\medskip

\section{Motivation and context }

On a filtred probability space $(\Omega ,\mathcal F,(\mathcal F_t)_{[0,1]},\mathbb P)$, we consider $X_1=(X_{1,t})_{t\in [0,1]}$ and $X_2=(X_{2,t})_{t\in [0,1]}$ two real processes defined by a L\'evy jump-diffusion constructed via the superposition of a Wiener process with drift and an independent compound Poisson process. This is one of the first and simplest extensions to the classical geometric Brownian motion underlying the famous Black-Scholes-Merton framework for option pricing.
\vspace{5pt}

More precisely, $X_1=(X_{1,t})_{t\in [0,1]}$ and $X_2=(X_{2,t})_{t\in [0,1]}$ are given by
\begin{equation}\label{equa1}
\left\{\aligned
dX_{1,t}=b_{1} (t,\omega)dt +\sigma_{1,t}dW_{1,t}+dJ_{1,t}\\
dX_{2,t}=b_{2}(t,\omega) dt+\sigma_{2,t} dW_{2,t}+dJ_{2,t}
\endaligned
\right.
\end{equation}
for $t\in [0,1]$ where  $W_1=(W_{1,t})_{t\in [0,1]}$ and $W_2=(W_{2,t})_{t\in [0,1]}$ are two correlated Wiener processes, with $\rho_t={\rm Cov}(W_{1,t},W_{2,t}), t\in [0,1]$. We can write $W_{2,t}=\rho_tdW_{1,t}+\sqrt{1-\rho_t^2}dW_{3,t},$
where $W_1=(W_{1,t})_{t\in [0,1]}$ and  $W_3=(W_{3,t})_{t\in [0,1]}$ are independent Wiener processes. $J_1$ and $J_2$ are possibly correlated pure jump processes. We assume here that $J_1$ and $J_2$  have  finite jump activity, that is a.s. there are only finitely many jumps on any  finite time interval. A general L\'evy model would contain also a compensated infinte activity pure jump component.
\vspace{5pt}

Under our assumption $J_{\ell}$ is necessarily a compound Poisson processe and it can be written as
$$J_{\ell,s}=\sum_{i=1}^{N_{\ell,s}}Y_{\ell,i},\qquad s\in [0,1].$$
Here $Y_{\ell,i}$ are i.i.d. real random variables having law $\nu_{\ell}/\lambda_{\ell}$, where $\nu_{\ell}$ is the L\'evy measure of $X_{\ell}$ normalized by the total mass $\lambda_{\ell}=\nu_{\ell}(\mathbb R-\{0\})<+\infty$, and $N_{\ell}$ is a poisson process, independent of each $Y_{\ell,i}$, and with constant intensity $\lambda_{\ell}$.

\vspace{10pt}
Such a jump-type stochastic process is recently a standard tool, e.g., for modeling asset values in finance and insurance. The key motivation behind jump-diffusion models is the incorporation of market "stocks", which result in "large" and sudden changes in the price of risky security and which can hardly be modeled by the diffusive component.

\vspace{10pt}
In this paper we concentrate on the estimation of
$$[\mathcal V]_t=\left(\int_0^t\sigma^2_{1,s}ds,\int_0^t\sigma^2_{2,s}ds,\int_0^t\sigma_{1,s}\sigma_{2,s}\rho_sds\right)$$

\vspace{5pt}
Over the last decade, several estimation methods for the integrated variance-covariance $\mathcal V_t$ have been proposed. We adopt the threshold estimator which is introduced  by Mancini \cite{Mancini4} and also by Shimizu and Yoshida \cite{shi3}, independently.

\vspace{5pt}
In this method, only the variations under a given threshold function are taken into account. The specific estimator excludes all terms containing jumps from the realized co-variation while remaining consistent, efficient and robust when synchronous data are considered.

\vspace{5pt}
Since the seminal work of Mancini  \cite{Mancini4}, several authors have leveraged or extended the thresholding cencept to deal with complex stochastic models, see Shimizu and Yoshida \cite{shi3}, or Ogihara and Yoshida \cite{OY}. The similar idea is also used by various authors in different contexts; see, e.g., A\"\i t-Sahalia et al. \cite{Sahalia4}, \cite{Sahalia2} and \cite{Sahalia3}, Gobbi and Mancini \cite{Mancini3} , Cont and Mancini \cite{Mancini2} , among others.

\vspace{10pt}

So, given the synchronous and evenly-spaced observation of the process $X_{1,t_0},X_{1,t_1},\cdots, X_{1,t_n},$ $X_{2,t_0},X_{2,t_1}\cdots,X_{2,t_n}$ with $t_0=0, t_n=1,n\in \mathbb N$, we consider  the following statistics
$$
\left(\sum_{k=1}^{[nt]}(\Delta^n_kX_{1})^2, \sum_{k=1}^{[nt]}(\Delta^n_kX_{2})^2, \sum_{k=1}^{[nt]}\Delta^n_kX_{1}\Delta^n_kX_{2} \right)
$$
where $\Delta^n_kX_{\ell}:=X_{\ell,t_k}-X_{\ell,t_{k-1}}$. However this estimate can be highly biased when the processes $X_{\ell}$ contain jumps, in fact, as $n\rightarrow \infty$ such a sum approaches the global quadratic variance-covariation
$$\left([X_1]_t,[X_2]_t,[X_1,X_2]_t\right)$$
where
$$[X_{\ell}]_t:=\int_0^t\sigma_{\ell,s}^2ds+\sum_{s\le t}(\Delta J_{\ell,s})^2,\quad {\rm and} \quad [X_1,X_2]_t:=\int_0^t\sigma_{1,s}\sigma_{2,s}\rho_sds+\sum_{s\le t}\Delta J_{1,s} \Delta J_{2,s}.$$
which also contain the co-jumps, where $\Delta J_{\ell,s}=J_{\ell,s}-J_{\ell,s^-}$.

If we take a deterministic function $r(\frac 1n)$ at the step $\frac 1n$ between the observations, such that
$$\lim_{n\rightarrow \infty}r\left(\dfrac 1n\right)=0, \qquad {\rm and}\quad \lim_{n\rightarrow \infty}\frac{\log n}{nr\left(\frac 1n\right)}=0.
$$
The function $r(\cdot)$ is a threshold such that whenever $|\Delta^n_kX_{\ell}|^2>r(\frac 1n)$,
a jump has to occur within $]t_{k-1},t_k]$. Hence we can recover $[\mathcal V]_t$ using the following threshold estimator

$$
\mathcal V_t^n(X)=(\mathcal Q_{1,t}^n(X), \mathcal Q_{2,t}^n(X),\mathcal C_t^n(X))
$$
where
$$
\mathcal Q_{\ell,t}^n(X)=\sum_{k=1}^{[nt]}(\Delta^n_kX_{\ell})^2{\bf 1}_{\{(\Delta^n_kX_{\ell})^2\le r(\frac 1n)\}}
$$
and
$$
\mathcal C_t^n(X)= \sum_{k=1}^{[nt]}\Delta^n_kX_{1} \Delta^n_kX_{2}{\bf 1}_{\{\max_{\ell=1}^2(\Delta^n_kX_{\ell})^2\le r(\frac 1n)\}}
$$

In the work \cite{JFLN}, the authors determine what constitutes a good threshold sequence $r_n$ and they propose an objective method for selecting such a sequence.

\vspace{5pt}
In the case that $X_{\ell}$ have no jumps, this question has been well investigated. The problem of the large deviation of the quadratic estimator of the integrated volatility (without jumps and in the case of synchronous sampling scheme) is obtained in the paper by Djellout et al. \cite{Djellout1} and recently Djellout and Samoura \cite{Djellout2} have studied the large deviation for the covariance estimator. Djellout et al. \cite{Djellout3} have also investigated the problem of the large deviation for the realized (co-)volatility vector which allows them to provide the large deviation for the standard dependence measures between the two assets returns such as the realized regression coefficients, or the realized correlation.
\vspace{5pt}

However, the inclusion of jumps within financial models seems to be more and more necessary for pratical applications. In this case,  Mancini \cite{Mancini2} has shown that $\mathcal V_t^n$ is a consistent estimators of  $\mathcal V_t$ and has some asymtotic normality respectively. Furthermore, when $\sigma_t=\sigma$, she \cite{Mancini1} studied the large deviation for the threshold estimator. Jiang \cite{Jiang1} obtained moderate deviations and functional moderate deviations for threshold estimator. In our paper and by the method as in Mancini \cite{Mancini1} and Djellout et al \cite{Djellout3}, we consider moderate and functionnal moderate deviation for estimators $V_t^n$ and large deviation.
\vspace{5pt}

More precisely we are interested in the estimations of
$$\mathbb P\left(\frac{\sqrt n}{v_n}\left(\mathcal V_t^n(X)-[\mathcal V]_t\right)\in A\right)$$
where $A$ is a given domain of deviation, $(v_n)_{n>0}$ is some sequence denoting the scale of deviation. When $v_n=1$ this is exactly the estimation of central limit theorem. When $v_n=\sqrt n$, it becomes the large deviation. Furthermore, when $v_n\rightarrow \infty$ and $v_n=o(\sqrt n)$, this is the so called moderate deviations. In other words, the moderate deviations investigate the convergence speed between the large deviations and central limit theorem.

\vspace{5pt}

Let us recall some basic defintions in large deviations theory. Let $(\mu_t)_{t>0}$ be a family of probability on a topological space $(S, {\mathcal S})$  where $\mathcal S$ is a $\sigma$-algebra on $S$ and $\lambda_t$ be a nonnegative function on $[1,+\infty[$ such that $\lim_{t\rightarrow \infty}\lambda_t=+\infty.$  A function $I:S\rightarrow [0,+\infty]$ is said to be a rate function if it is lower semicontinuous and it is said to be a good rate function if its level set $\{x\in S ; I(x)\le a\}$ is a compact for all $a\ge 0$.

$(\mu_t)$ is said to satisfy a large deviation principle with speed $\lambda_t$ and rate function $I$ if for any closed set $F\in \mathcal S$
$$\limsup_{t\rightarrow \infty }\frac{1}{\lambda_t}\log\mu_t(F)\le -\inf_{x\in F}I(x)$$
and for any open set $G\in \mathcal S$
$$\limsup_{t\rightarrow \infty }\frac{1}{\lambda_t}\log\mu_t(G)\ge -\inf_{x\in G}I(x).$$

\noindent{\bf Notations.} \textit{In the whole paper, for any matrix $M$, $M^T$ and $\Vert M \Vert$ stand for the transpose and the euclidean norm of $M$, respectively. For any square matrix $M$, $\det(M)$ is the determinant of $M$. Moreover, we will shorten \textnormal{large deviation principle} by LDP and \textnormal{moderate deviation principle} by MDP.
We denote by $\langle\cdot,\cdot \rangle$ the usual scalar product. For any process $Z_t$, $\Delta_s^tZ$ stands for the increment $Z_t-Z_s$. We use $\Delta_k^nZ$ for $\Delta_{t_{k-1}^n}^{t_k^n}Z$. In addition, for a sequence of random variables $(Z_{n})_{n}$ on $\dR^{d \times p}$, we say that $(Z_{n})_{n}$ converges $(\lambda_{n})-$superexponentially fast in probability to some random variable $Z$ if, for all $\delta > 0$,
\begin{equation*}
\limsup_{n \rightarrow \infty} \frac{1}{\lambda_n} \log \dP\Big( \left\Vert Z_{n} - Z \right\Vert > \delta \Big) = -\infty.
\end{equation*}
This \textnormal{exponential convergence} with speed $\lambda_{n}$ will be shortened as
\begin{equation*}
Z_{n} \superexp Z.
\end{equation*}}

\vspace{10pt}
The article is arranged in two upcoming sections. Section 2 is devoted to our main results on the LDP and MDP for the (co-)volatility vector in the presence of jumps. In section 3, we give the proof of these theorems.

\section{Main results }

Let $X_t=(X_{1,t},X_{2,t})$ be given by \eqref{equa1}. We introduce the following conditions

\vspace{10pt}
{\bf (B)} for $\ell=1,2$ $b(\cdot,\cdot)\in L^\infty(dt\otimes \mathbb{P})$

\vspace{10pt}
{\bf (LDP)} Assume that for $\ell=1,2$
\begin{itemize}
\item $\sigma_{\ell,t}^2(1-\rho_t^2)$ and $\sigma_{1,t}\sigma_{2,t}(1-\rho_t^2)$ $\in$ $L^\infty([0,1],dt)$.
\item the functions $t \to \sigma_{\ell,t}$ and $t \to \rho_t$ are continuous.
\item let $r$ such that
$$r\left(\dfrac 1n\right)\xrightarrow[{n \to \infty}]\quad 0 \quad {\rm and} \quad nr\left(\dfrac 1n\right)\xrightarrow[{n \to \infty}]\quad \infty.$$
\end{itemize}

\vspace{10pt}
{\bf (MDP)} Assume that for $\ell=1,2$
\begin{itemize}
\item $\sigma_{\ell,t}^2(1-\rho_t^2)$ and $\sigma_{1,t}\sigma_{2,t}(1-\rho_t^2)$ $\in$ $L^2([0,1],dt)$.
\item Let $(v_n)_{n\geqslant1}$ be a sequence of positive
numbers such that
\begin{eqnarray}\label{Max_condition}
&&v_n \xrightarrow[{n \to \infty}]\quad  \infty \quad {\rm and} \quad \dfrac{v_n}{\sqrt{n}} \xrightarrow[{n \to \infty}]\quad 0 \quad {\rm and} \quad \sqrt n v_nr\left(\frac 1n\right)=O(1)\nonumber \\
&&{\rm and \,\,for}\quad \ell=1,2\quad\frac{r\left(\dfrac 1n\right)}{\log\left(\dfrac{n}{v_n^2}\right)\displaystyle\max_{k=1}^n\int_{t_{k-1}}^{t_k}\sigma_{\ell,s}^2ds}\longrightarrow +\infty.
\end{eqnarray}
\end{itemize}

\vspace{10pt}

We introduce the following function, which will play a crucial role in the calculation of the moment generating function: for $-1<c<1$ let for any $\lambda=(\lambda_1,\lambda_2,\lambda_3)\in \mathbb{R}^3$
\begin{equation}\label{ASSA}
P_c(\lambda):=
  \left\{
  \begin{array}{lcl}
  \vspace{0.3cm}
 -\dfrac{1}{2} \log\left(\dfrac{(1-2\lambda_1(1-c^2))(1-2\lambda_2(1-c^2))-(\lambda_3(1-c^2)+c)^2}{1-c^2}\right)\\
 \vspace{0.3cm}
\qquad \qquad\qquad \qquad \qquad\qquad if\qquad\lambda\in {\mathcal D}\\

  +\infty, \quad otherwise
  \end{array}
  \right.
\end{equation}
where
\begin{equation}\label{defD}
\displaystyle {\mathcal D}_c=\left\{\lambda\in \mathbb{R}^3,\,\,\max_{\ell=1,2}\lambda_{\ell} < \dfrac{1}{2(1-c^2)}\,\,{\rm and}\,\prod_{\ell =1}^2\left(1-2\lambda_{\ell}(1-c^2)\right)> \left(\lambda_3(1-c^2)+c\right)^2\right\}.
\end{equation}
\vspace{10pt}

Let us present now the main results.

\subsection{Moderate deviation}

Let us now consider the intermediate scale between the central limit theorem and the law of large numbers.
\begin{thm}\label{MDP1} For t=1 fixed. Under the conditions {\bf (MDP)} and {\bf (B)}, the sequence
 $$\dfrac{\sqrt{n}}{v_n}\left(\mathcal V_{1}^n(X)-[\mathcal V]_{1}\right)$$
 satisfies the LDP on $\mathbb{R}^3$ with speed $v_n^2$ and with rate function given by
 \begin{equation}\label{MDPtaux1}
 I_{mdp}(x)=\sup_{\lambda\in\mathbb{R}^3}\left(\left\langle \lambda,x \right\rangle - \dfrac{1}{2}\left\langle\lambda, \Sigma_1\cdot\lambda\right\rangle\right)=\frac{1}{2}\left\langle x, \Sigma_1^{-1}\cdot x\right\rangle
 \end{equation}
 with
 $$
 \Sigma_{1}= \begin{pmatrix}
                                       \vspace{0.3cm}
                                       \int_0^1\sigma_{1,t}^4 \mathrm dt &  \int_0^1\sigma_{1,t}^2\sigma_{2,t}^2\rho_t^2 \mathrm dt & \int_0^1\sigma_{1,t}^3\sigma_{2,t}\rho_t \mathrm dt\\
                                        \vspace{0.3cm}
                                        \int_0^1\sigma_{1,t}^2\sigma_{2,t}^2\rho_t^2 \mathrm dt &  \int_0^1\sigma_{2,t}^4 \mathrm dt & \int_0^1\sigma_{1,t}\sigma_{2,t}^3\rho_t \mathrm dt\\
                                        \int_0^1\sigma_{1,t}^3\sigma_{2,t}\rho_t \mathrm dt & \int_0^1\sigma_{1,t}\sigma_{2,t}^3\rho_t \mathrm dt &  \int_0^1\dfrac{1}{2}\sigma_{1,t}^2\sigma_{2,t}^2(1+\rho_t^2) \mathrm dt
                                      \end{pmatrix}.
 $$
\end{thm}

\vspace{10pt}
\begin{rem} Under the condition $b_{\ell}=0$, we can prove that  for all $\theta \in \mathbb R^3$
$$
\lim_{n\to\infty}\frac{1}{v_n^2}\log\mathbb E\left(e^{\sqrt{n}v_n\langle\theta,\mathcal V_1^n(X)-[\mathcal V]_1\rangle}\right)=\frac{1}{2}<\theta, \Sigma_1\cdot \theta>.
$$
This gives an alternative proof of the moderate deviation using G\"artner-Ellis theorem.
\end{rem}

 \begin{rem}
 If for some $p>2$, $\sigma_{1,t}^2$, $\sigma_{2,t}^2$ and $\sigma_{1,t}\sigma_{2,t}(1-\rho_t^2)$ $\in$ $L^p([0,1])$ and $v_n=O(n^{\frac{1}{2}-\frac{1}{p}})$,
  the condition \eqref{Max_condition} in {\bf (MDP)} is verified.
 \end{rem}
\vspace{10pt}

 Let ${\mathcal H}$ be the banach space of $\mathbb{R}^3$-valued  right-continuous-left-limit non decreasing functions $\gamma$ on $[0,1]$ with $\gamma(0)=0$, equipped with the uniform norm
and the $\sigma-$field $\mathcal{B}^s$ generated by the coordinate $\{\gamma(t),0 \leqslant t \leqslant 1\}$.

\begin{thm}
 \label{MDP2}
 Under the conditions {\bf (MDP)} and {\bf (B)}, the sequence
 $$\dfrac{\sqrt{n}}{v_n}\left(\mathcal V_{.}^n(X)-[\mathcal V]_{.}\right)$$
 satisfies the LDP on $\mathcal H$ with speed $v_n^2$ and with rate function given by
 \begin{equation}
 J_{mdp}(\phi)=\label{ARR}
  \left\{
  \begin{array}{lcl}
  \vspace{0.3cm}
 \displaystyle \int_0^1 \frac{1}{2}\left\langle \dot \phi(t),\overline{\Sigma}_t^{-1}\cdot \dot \phi(t) \right\rangle dt\qquad if \quad \phi \in \mathcal{AC}_0([0,1])\\		
  \vspace{0.3cm}
  +\infty, \qquad otherwise,
  \end{array}
  \right.
  \end{equation}
 where
 $$
 \overline{\Sigma}_t= \begin{pmatrix}
                                       \vspace{0.3cm}
                                       \sigma_{1,t}^4 & \sigma_{1,t}^2\sigma_{2,t}^2\rho_t^2 & \sigma_{1,t}^3\sigma_{2,t}\rho_t \\
                                        \vspace{0.3cm}
                                        \sigma_{1,t}^2\sigma_{2,t}^2\rho_t^2 & \sigma_{2,t}^4 & \sigma_{1,t}\sigma_{2,t}^3\rho_t \\
                                        \sigma_{1,t}^3\sigma_{2,t}\rho_t & \sigma_{1,t}\sigma_{2,t}^3\rho_t & \dfrac{1}{2}\sigma_{1,t}^2\sigma_{2,t}^2(1+\rho_t^2)
                                      \end{pmatrix}
 $$
 is invertible and  $\overline{\Sigma}_t^{-1}$ his inverse such that
 $$
 \overline{\Sigma}_t^{-1}=\frac{1}{{\rm det}(\overline{\Sigma}_t)} \begin{pmatrix}
												  \vspace{0.3cm}
												  \dfrac{1}{2}\sigma_{1,t}^2\sigma_{2,t}^6(1-\rho_t^2) & \dfrac{1}{2}\sigma_{1,t}^4\sigma_{2,t}^4\rho_t^2(1-\rho_t^2) & -\sigma_{1,t}^3\sigma_{2,t}^5\rho_t(1-\rho_t^2)\\
												    \vspace{0.3cm}
												    \dfrac{1}{2}\sigma_{1,t}^4\sigma_{2,t}^4\rho_t^2(1-\rho_t^2) & \dfrac{1}{2}\sigma_{1,t}^6\sigma_{2,t}^2(1-\rho_t^2) & -\sigma_{1,t}^5\sigma_{2,t}^3\rho_t(1-\rho_t^2)\\
												    -\sigma_{1,t}^3\sigma_{2,t}^5\rho_t(1-\rho_t^2) & -\sigma_{1,t}^5\sigma_{2,t}^3\rho_t(1-\rho_t^2) &  \sigma_{1,t}^4\sigma_{2,t}^4(1-\rho_t^4)
												  \end{pmatrix},
$$												

$$with\qquad {\rm det}(\overline{\Sigma}_t)=\dfrac{1}{2}\sigma_{1,t}^6\sigma_{2,t}^6(1-\rho_t^2)^3,$$
and $\mathcal{AC}_0=\left\{ \phi:[0,1] \rightarrow \mathbb{R}^3\,\, {\it  is\, absolutely\, continuous\, with}\,\phi(0)=0\right\}.$
\end{thm}
\vspace{10pt}

\begin{rem} A similar result for the moderate deviations is obtained by Jiang \cite{Jiang1}  in the jump case for $\left(\frac{\sqrt n}{v_n}\left(\mathcal Q_{\ell,t}^n-\int_0^t\sigma^2_{\ell,s}ds\right)\right)_{n\ge 1}.$
\end{rem}

\vspace{10pt}

\subsection{Large deviation}

Our second result is about the large deviation of $\mathcal V_{1}^n (X)$, i.e. at fixed time.

\begin{thm}\label{LDP1} Let $t=1$ be fixed. Under the conditions {\bf (LDP)} and {\bf (B)}  , the sequence $\mathcal V_{1}^n(X)$ satisfies the LDP on $\mathbb{R}^3$ with speed $n$ and with good rate function given by the legendre transformation of $\Lambda$, that is
  \begin{equation}\label{LDPtaux1}
  I_{ldp}(x)=\sup_{\lambda\in\mathbb{R}^3}\left(\left\langle \lambda,x \right\rangle - \Lambda(\lambda)\right),\end{equation}
where~$\Lambda(\lambda)=\int_0^1P_{\rho_t}(\lambda_1\sigma_{1,t}^2,\lambda_2\sigma_{2,t}^2,\lambda_3\sigma_{1,t}\sigma_{2,t})dt$.
\end{thm}
\vspace{10pt}

\begin{rem} Under the condition $b_{\ell}=0$, we can calculate the moment generating function of $ \mathcal V_1^n(X)$. We obtain that for all  $\theta=(\theta_1,\theta_2,\theta_3)^T\in {\mathcal D}_{\rho_t}$
$$
\lim_{n\to\infty}\frac{1}{n}\mathbb E\left(e^{n\langle\theta, \mathcal V_1^n(X)\rangle}\right)
=\int_0^1P_{\rho_s}\left(\theta_1\sigma_{1,s}^{2},\theta_2\sigma_{2,s}^{2},\theta_3\sigma_{1,s}\sigma_{2,s}\right)ds.
$$
But the study of the steepness is more difficult.
\end{rem}

Let us consider the case where diffusion and correlation coefficients are constant, the rate function being easier to read. Before that let us introduce the function $P^*_c$ which is the Legendre transformation of $P_c$ given in (\ref{ASSA}), for all $x=(x_1,x_2,x_3)$
\begin{equation}\label{Legendre}
P^*_c(x):=
  \left\{
  \begin{array}{lcl}
  \vspace{0.3cm}
  \log\left(\dfrac{\sqrt{1-c^2}}{\sqrt{x_1x_2-x_3^2}}\right)-1+\dfrac{x_1+x_2-2cx_3}{2(1-c^2)}\\
 \vspace{0.3cm}
\qquad \qquad\qquad \qquad if\qquad x_1>0,\,\,x_2>0,\,\,x_1x_2>x_3^2\\

  +\infty, \quad otherwise.
  \end{array}
  \right.
\end{equation}

\begin{cor} \label{LDPconstant}
We assume that  for $\ell=1,2$ $\sigma_{\ell}$ and $\rho$ are constants. Under the condition {\bf (B)},  we obtain that $\mathcal V_{1}^n(X)$ satisfies the LDP on $\mathbb{R}^3$ with speed $n$ and with good rate function $I_{ldp}^{\mathcal V}$ given by
\begin{equation}\label{LDPtauxconstantV}
I_{ldp}^{\mathcal V}(x_1,x_2,x_3)=P^*_{\rho}\left(\frac{x_1}{\sigma_1^2},\frac{x_2}{\sigma_2^2}, \frac{x_3}{\sigma_1\sigma_2 }\right),
\end{equation}
where $P^*_c$ is given in (\ref{Legendre}).
\end{cor}

\begin{rem} In the case $\sigma_{\ell}$ is constant, a similar result for the large deviations is obtained by Mancini \cite{Mancini1}  in the jump case for $\left(\mathcal Q_{\ell,1}^n\right)_{n\ge 1}$
\end{rem}
\vspace{10pt}

Now, we shall extend Theorem \ref{LDP1} to the process-level large deviations, i.e. for trajectories
$(\mathcal V_{t}^n(X), t\in[0,1])$ which is interesting from the viewpoint of non-parametric statistics.

\vspace{10pt}

Let $\mathcal BV([0,1],\mathbb{R}^3)$~(shorted in $\mathcal BV$) be the space of functions of bounded variation on
$[0,1]$. We identify $\mathcal BV$ with $\mathcal M_3([0,1])$, the set of vector measures with value in $\mathbb{R}^3$.
This is done in the usual manner: to $f\in\mathcal BV$, there corresponds~$\mu^f$ by $\mu^f([0,t])=f(t)$.
Up to this identification, $\mathcal C_3([0,1])$ the set of $\mathbb{R}^3$-valued continuous bounded functions on $[0,1]$,
is the topology dual of $\mathcal BV$. We endow $\mathcal BV$ with the weak-* convergence topology
$\sigma\left(\mathcal BV, \mathcal C_3([0,1])\right)$ and with the associated Borel-$\sigma$-field $\mathcal B_{\omega}$.
Let $f\in\mathcal BV$ and $\mu^f$ the associated measure in $\mathcal M_3([0,1])$. Consider the Lebesgue decomposition of
$\mu^f$, $\mu^f=\mu^f_a+\mu^f_s$ where $\mu^f_a$ denotes the absolutely continuous part of $\mu^f$ with respect to $dx$ and
$\mu^f_s$ its singular part. We denote by $f_a(t)=\mu^f_a([0,t])$ and by $f_s(t)=\mu_s^f([0,t])$.

\begin{thm}\label{LDP2}
Under the conditions {\bf (LDP)} and {\bf (B)}, the sequence
 $\mathcal V_{.}^n(X)$
 satisfies the LDP on $\mathcal BV$ with speed $n$ and rate function $J_{ldp}$ given for any $f=(f_1,f_2,f_3)\in\mathcal BV$ by
 \begin{align}\label{LDP2taux}
 J_{ldp}(f)&=\int_0^1 P^*_{\rho_t}\left(\frac{f'_{1,a}(t)}{\sigma^2_{1,t}},\frac{f'_{2,a}(t)}{\sigma^2_{2,t}},
 \frac{f'_{3,a}(t)}{\sigma_{1,t}\sigma_{2,t}}\right)\\
 &+\int_0^1\frac{\sigma^2_{2,t}f'_{1,s}(t)+\sigma^2_{1,t}f'_{2,s}(t)-2\rho_t\sigma_{1,t}\sigma_{2,t}f'_{3,s}(t)}{2\sigma^2_{1,t}\sigma^2_{2,t}(1-\rho_t^2)}
 {\bf 1}_{[t:f'_{1,s}>0,f'_{2,s}>0,(f'_{3,s})^2<f'_{1,s}f'_{2,s}]}d\theta(t),\nonumber
 \end{align}
 where $P^*_c$ is given in (\ref{Legendre}) and $\theta$ is any real-valued nonnegative measure with respect
 to which $\mu_s^f$ is absolutely continuous and $f'_s={d\mu^f_s}/{d\theta}=(f'_{1,s},f'_{2,s},f'_{3,s})$.
\end{thm}

\vspace{10pt}

\section{Proofs }

For the convenience of the reader, we recall the following lemma which is the key of the proofs.
\begin{lem}\label{Approximation Lemma} (Approximation Lemma) Theorem 4.2.13 in  \cite{Demzei}

Let $(Y^n, X^n,n\in \mathbb N)$ be a family of random varibales valued in a Polish space $S$ with metric $d(\cdot,\cdot)$, defined on a probability space $(\Omega, {\mathcal F}, \mathbb P)$. Assume
\begin{itemize}
\item $\mathbb P(Y^n\in \cdot)$ satisfies the large deviation principle with speed $\epsilon_n$ ($\epsilon_n\rightarrow \infty$) and the good rate function $I$.
\item for every $\delta>0$
$$\limsup_{n\rightarrow \infty}\frac{1}{\epsilon_n}\log\mathbb P(d(Y^n,X^n)>\delta)=-\infty.$$
\end{itemize}
Then $\mathbb P(X^n\in \cdot)$ satisfies the large deviation principle with speed $\epsilon_n$ and the good rate function $I$.
\end{lem}

Before starting the proof, we need to introduce some technical tools. In the case without jumps, we introduce the following diffusion for $\ell=1,2$
$$D_{\ell,t} =\int_0^t\sigma_{\ell,s}dW_{\ell,s},$$
where $W_{\ell,s}$ and $\sigma_{\ell,s}$ are defined as before. We introduce the correspondent estimator
$$V_t^n=(Q_{1,t}^n,Q_{2,t}^n, C_t^n)$$
where for $\ell=1,2$
$$Q_{\ell,t}^n=\sum_{k=1}^{[nt]}\left(\Delta^n_kD_{\ell}\right)^2\quad {\rm and}\quad C_t^n=\sum_{k=1}^{[nt]}\Delta^n_kD_{1}\Delta^n_kD_{2}.$$

We recall the following results from Djellout et al. \cite{Djellout3}

\begin{prop}
Under the conditions  {\bf (B)} and {\bf (MDP)},

\begin{enumerate}
\item the sequence
 $$\dfrac{\sqrt{n}}{v_n}\left(V_{1}^n-[\mathcal V]_{1}\right)$$
 satisfies the LDP on $\mathbb{R}^3$ with speed $v_n^2$ and with rate function given by (\ref{MDP1}).

\item the sequence
 $$\dfrac{\sqrt{n}}{v_n}\left( V_{\cdot}^n-[\mathcal V]_{\cdot}\right)$$
 satisfies the LDP on $\mathcal H$ with speed $v_n^2$ and with rate function given by (\ref{MDP2}).
 \end{enumerate}
 \end{prop}

 \begin{prop}
 Under the conditions {\bf (B)} and {\bf (LDP)},
 \begin{enumerate}

\item the sequence $V_{1}^n$ satisfies the LDP on $\mathbb{R}^3$ with speed $n$ and with good rate function given in (\ref{LDPtaux1}).

\item the sequence
 $\mathcal V_{.}^n$
 satisfies the LDP on $\mathcal BV$ with speed $n$ and rate function $J_{ldp}$ given by (\ref{LDP2taux}).
 \end{enumerate}
\end{prop}

\vspace{10pt}


\subsection{Proof  of Theorem \ref{MDP1}}$\qquad$
\vspace{10pt}

We will do the proof in two steps.
\vspace{10pt}

{\underline {\it Part 1}} We start with the case $b_{\ell}=0$. In this case, $\mathcal V_t^n(X)=\mathcal V_t^n(X^0)$
with $X^0_{\ell,t}=X_{\ell,t}-\int_0^tb_{\ell}(s,\omega)ds$ and
$$
\mathcal Q_{\ell,1}^n(X^0)=\sum_{k=1}^{n}(\Delta^n_kX^0_{\ell})^2
{\bf 1}_{\{(\Delta^n_kX^0_{\ell})^2\le r(\frac 1n)\}}, ~~\ell=1,2
$$
and
$$
\mathcal C_1^n(X^0)= \sum_{k=1}^{n}\Delta^n_kX_{1}^0
\Delta^n_kX_{2}^0{\bf 1}_{\{\max_{\ell =1}^2(\Delta^n_kX_{\ell}^0)^2\le r(\frac 1n)\}}.
$$

We will prove that
$$\frac{\sqrt n}{v_n}\left(\mathcal V_1^n(X^0)-V_1^n\right)\superexpmdp0.$$

For that, we will prove  that for $\ell =1,2$
\begin{equation}\label{nq}
\frac{\sqrt n}{v_n}\left(\mathcal Q_{\ell,1}^n(X^0)-Q_{\ell,1}^n \right)\superexpmdp 0,
\end{equation}
and
\begin{equation}\label{nc}
\frac{\sqrt n}{v_n}\left(\mathcal C_1^n(X^0)- C_1^n\right)\superexpmdp 0.
\end{equation}

We start by the proof of (\ref{nq}). Since the processes $X^0_{\ell}$ and $D_{\ell}$ have independent increment, by Chebyshev inequality we obtain for all $\theta>0$
$$
\mathbb P\left(\frac{\sqrt n}{v_n}\left( \mathcal Q_{\ell,1}^n(X^0)-Q_{\ell,1}^n \right)>\delta\right)\le e^{-\theta\delta v_n^2}\prod_{k=1}^{n}\mathbb E\left(e^{\theta\sqrt n v_n\left[(\Delta^n_kX^0_{\ell})^2{\bf 1}_{\{(\Delta^n_kX_{\ell})^2\le r(\frac 1n)\}}-(\Delta^n_kD_{\ell})^2\right]}\right).
$$

We have to control each term appearing in the product
\begin{eqnarray}\label{Q}
\mathbb E\left(e^{\theta\sqrt n v_n\left[(\Delta^n_kX^0_{\ell})^2{\bf 1}_{\{(\Delta^n_kX^0_{\ell})^2\le r(\frac 1n)\}}-(\Delta^n_kD_{\ell})^2\right]}\right) \leq\Re_1(k,n)+\Re_2(k,n),
\end{eqnarray}
where
$$\Re_1(k,n):= \mathbb E\left(e^{\theta\sqrt n v_n\left[(\Delta^n_kX^0_{\ell})^2-(\Delta^n_kD_{\ell})^2\right]}{\bf 1}_{\{(\Delta^n_kX^0_{\ell})^2\le r(\frac 1n)\}}\right)$$
and
 $$\Re_2(k,n):=\mathbb P\left((\Delta^n_kX^0_{\ell})^2> r(\frac 1n)\right).$$

For the first term, we write
\begin{eqnarray}\label{r0}
\Re_1(k,n)&=& \mathbb E\left(e^{\theta\sqrt n v_n\left[(\Delta^n_kX^0_{\ell})^2-(\Delta^n_kD_{\ell})^2\right]}{\bf 1}_{\{(\Delta^n_kX^0_{\ell})^2\le r(\frac 1n)\}}|\Delta^n_kN_{\ell}=0\right) \mathbb P(\Delta^n_kN_{\ell}=0)  \nonumber  \\
& &\quad +\mathbb E\left(e^{\theta\sqrt n v_n\left[(\Delta^n_kX^0_{\ell})^2-(\Delta^n_kD_{\ell})^2\right]}{\bf 1}_{\{(\Delta^n_kX^0_{\ell})^2\le r(\frac 1n), \Delta^n_kN_{\ell}\not=0\}}\right).
\end{eqnarray}

Since $N_{\ell}$ is independent of $W_{\ell}$, we obtain that
\begin{eqnarray}{\label{r1}}
\Re_1(k,n)&\le& \mathbb P\left((\Delta^n_kD_{\ell})^2\le r(\frac 1n)\right) e^{-\lambda_{\ell}/n} +e^{\sqrt n v_n \theta r(\frac 1n)} (1-e^{-\lambda_{\ell}/n})\nonumber\\
&\le&1+e^{\sqrt n v_n \theta r(\frac 1n)} (1-e^{-\lambda_{\ell}/n}).
\end{eqnarray}

Now we have to control $\Re_2(k,n)$, by the same argument as before we have
\begin{eqnarray*}
\Re_2(k,n) &=& \mathbb P\left((\Delta^n_kX^0_{\ell})^2> r(\frac 1n)|\Delta^n_kN_{\ell}=0\right) \mathbb P(\Delta^n_kN_{\ell}=0)   \nonumber \\
& &\quad +\mathbb P\left(\Delta^n_kX^0_{\ell})^2> r(\frac 1n), \Delta^n_kN_{\ell}\not=0\right)\nonumber \\
&\le& \mathbb P\left((\Delta^n_kD_{\ell})^2>r(\frac 1n) \right)e^{-\lambda_{\ell}/n}+ (1-e^{-\lambda_{\ell}/n}).
\end{eqnarray*}

From exponential inequality for martingales, it follows that for $\ell =1,2$,
\begin{equation}\label{expoential-ineq}
\mathbb P\left((\Delta^n_kD_{\ell})^2>r\left(\frac 1n\right) \right)\le \exp\left(-\frac{r(\frac 1n)}{2\int_{t_{k-1}}^{t_k}\sigma_{\ell,s}^2ds}\right),
\end{equation}
which implies that
\begin{eqnarray}{\label{r2}}
\Re_2(k,n) &\le&\exp\left(-\frac{r(\frac 1n)} {2\int_{t_{k-1}}^{t_k}\sigma_{\ell,s}^2ds}\right)+(1-e^{-\lambda_{\ell}/n}).
\end{eqnarray}

From (\ref{Q}), (\ref{r1}) and (\ref{r2}), we obtain that
\begin{eqnarray*}
\mathbb E\left(e^{\theta\sqrt n v_n\left[(\Delta^n_kX^0_{\ell})^2{\bf 1}_{\{(\Delta^n_kX^0_{\ell})^2\le r(\frac 1n)\}}-(\Delta^n_kD_{\ell})^2\right]}\right) &\le& 1+ ( 1+ e^{\sqrt n v_n \theta r(\frac 1n)}) (1-e^{-\lambda_{\ell}/n})\\
&&\qquad+\exp\left(-\frac{r(\frac 1n)}{2\int_{t_{k-1}}^{t_k}\sigma_{\ell,s}^2ds}\right).
\end{eqnarray*}

Using the hypotheses {\bf (MDP)}, we have
\begin{equation}\label{E-e-1}
\limsup_{n\rightarrow\infty}\frac{n}{v_n^2}
\max_{k=1}^n\log \mathbb E\left(e^{\theta\sqrt n v_n\left[(\Delta^n_kX^0_{\ell})^2{\bf 1}_{\{(\Delta^n_kX^0_{\ell})^2\le r(\frac 1n)\}}-(\Delta^n_kD_{\ell})^2\right]}\right)=0.
\end{equation}

So
$$
\limsup_{n\rightarrow\infty}\frac{1}{v_n^2}\log\mathbb P\left(\frac{\sqrt n}{v_n}
\left(\mathcal Q_{\ell,1}^n(X^0)-Q_{\ell,1}^n \right)>\delta\right)\le -\lambda \delta.
$$

Letting $\lambda$ goes to infinity, we obtain that the right hand of the last inequality goes to $-\infty$. Proceeding in the same way for $-(\mathcal Q_{\ell,t}^n(X^0)-Q_{\ell,t}^n)$ we obtain (\ref{nq}).

\vspace{10pt}

Now we have to prove (\ref{nc}). For that we have the following decompostion
 \begin{equation}\label{decompC}
\mathcal C_{1}^n(X^0)-C_{1}^n =\frac 12 \left[\overline{\mathcal Q}_{3,1}^n(X^0)- Q_{3,1}^n\right]-\frac 12\left[\sum_{\ell=1}^2\overline{\mathcal Q}_{\ell,1}^n(X^0)-Q_{\ell,1}^n\right],
\end{equation}
where
$$Q_{3,1}^n=\sum_{k=1}^n (\Delta^n_kD_{1}+\Delta^n_kD_{2})^2,$$
and for $\ell=1,2$
$$\overline{\mathcal Q}_{\ell,t}^n(X^0)=\sum_{k=1}^n(\Delta^n_kX^0_{\ell})^2 {\bf 1}_{\{\max_{\ell=1}^2(\Delta^n_kX^0_{\ell})^2\le r(\frac 1n)\}}$$
and
$$\overline{\mathcal Q}_{3,1}^n(X^0)=\sum_{k=1}^n(\Delta^n_kX^0_{1}+\Delta^n_kX^0_{2})^2 {\bf 1}_{\{\max_{\ell=1}^2(\Delta^n_kX^0_{\ell})^2\le r(\frac 1n)\}}.$$

Remark that $\overline{\mathcal Q}_{\ell,t}^n(X^0)$ is a slight modification of $\mathcal Q_{\ell,t}^n(X^0).$

We know that $\Delta^n_kD_{1}+\Delta^n_kD_{2}\sim {\mathcal N}(0,\beta^2(k,n))$ with $$
\beta^2(k,n)=\int_{t_{k-1}}^{t_k}\sigma_{1,s}^2ds+\int_{t_{k-1}}^{t_k}\sigma_{2,s}^2ds+2\int_{t_{k-1}}^{t_k}\sigma_{1,s}\sigma_{2,s}\rho_sds.
$$
For all $\delta>0$, we have
$$\mathbb P\left(\frac{\sqrt n}{v_n}\left|\mathcal C_{1}^n(X^0)-C_{1}^n \right|>\delta \right)\le 3\max_{\ell =1}^3\mathbb P \left(\frac{\sqrt n}{v_n}\left| \overline{\mathcal Q}_{\ell,1}^n(X^0)-Q_{\ell,1}^n \right|>\frac{2\delta}{3}\right).$$

So we obtain (\ref{nc}).


\vspace{10pt}

{\underline{\it Part 2}} We have to prove that
$$\frac{\sqrt n}{v_n}\left(\mathcal V_1^n(X)-\mathcal V_1^n(X^0)\right)\superexpmdp 0.$$

We have that
\begin{equation}\label{QQ}\left|\mathcal Q_{\ell,1}^n(X)-\mathcal Q_{\ell,1}^n(X^0)\right|\le \varepsilon(n) \mathcal Q_{\ell,1}^n(X^0)+\left(1+\frac{1}{\varepsilon(n)}\right)Z^n_{\ell}\end{equation}
and
\begin{equation}\label{C}\left|\mathcal C_{1}^n(X)-\mathcal C_{1}^n(X^0)\right|\le \varepsilon(n) \max_{\ell=1}^2\mathcal Q_{\ell,1}^n(X^0)+\left(1+\frac{1}{\varepsilon(n)}\right)\max_{\ell=1}^2Z^n_{\ell},\end{equation}
where
$$Z^n_{\ell}=\sum_{k=1}^n\left(\int_{t_{k-1}}^{t_k}b_{\ell}(s,\omega)ds\right)^2.$$
By the condition {\bf (B)}, we have that $\left\|Z_{\ell}^n\right\|\le \dfrac 1n.$ We choose $\varepsilon(n)$ such that 
$$
\frac{\sqrt {n}}{v_n}\varepsilon(n)\rightarrow 0,\quad v_n\sqrt{n}\varepsilon(n)\to\infty,
$$ 
so by the MDP of  $\mathcal Q_{\ell,1}^n(X^0)$, we obtain the result.

\vspace{10pt}

\subsection{Proof  of Theorem \ref{MDP2}}$\qquad$
\vspace{10pt}

Since the sequence $\frac{\sqrt{n}}{v_n}(V_{\cdot}^n-[\mathcal V]_{\cdot})$ satisfies the LDP on $\mathcal H$
with speed $v_n^2$ and rate function~$J_{mdp}$, by Lemma \ref{Approximation Lemma}, it is sufficient to show that:
\begin{equation}\label{exponential-convergence-mdp}
\frac{\sqrt{n}}{v_n}\sup_{t\in[0,1]}\left\|\mathcal{V}_{t}^n(X^0)-V_{t}^n\right\|\superexpmdp0.
\end{equation}

\begin{lem}\label{lem-mdp-1}
Under the condition {\bf (MDP)}, we have
$$
\lim_{n\to\infty}\frac{\sqrt{n}}{v_n}\sup_{t\in[0,1]}\left\|\mathbb E\mathcal V_t^n(X^0)-[\mathcal V]_t\right\|=0.
$$
\end{lem}
\noindent\textbf{Proof} We will prove that for $\ell=1,2$
\begin{equation}\label{e-1}
\lim_{n\to\infty}\frac{\sqrt{n}}{v_n}\sup_{t\in[0,1]}\left|\mathbb E\mathcal Q_{\ell,t}^n(X^0)-\int_0^t\sigma_{\ell,s}^2ds\right|=0.
\end{equation}
and
\begin{equation}\label{e-2}
\lim_{n\to\infty}\frac{\sqrt{n}}{v_n}\sup_{t\in[0,1]}\left|\mathbb E\mathcal C_{t}^n(X^0)-\int_0^t\sigma_{1,s}\sigma_{1,s}\rho_sds\right|=0.
\end{equation}

In fact, (\ref{e-1}) can be done in the same way as in Jiang \cite{Jiang1}. It remains to show (\ref{e-2}).
Using (\ref{decompC}), we obtain that

\begin{align*}
\left|\mathbb E\mathcal C_{t}^n(X^0)-\int_0^t\sigma_{1,s}\sigma_{1,s}\rho_sds\right|\le
\frac 12 \left|\mathbb E\overline{\mathcal Q}_{3,t}^n(X^0)- \beta_{t}\right|+\max_{\ell=1}^2\left|\mathbb E\overline{\mathcal Q}_{\ell,t}^n(X^0)-\int_0^t\sigma^2_{\ell,s}ds\right|,
\end{align*}
where $\beta_t=\int_0^t\sigma^2_{1,s}ds+\int_0^t\sigma^2_{2,s}ds+2\int_0^t\sigma_{1,s}\sigma_{2,s}\rho_sds.$ So the proof of  (\ref{e-2}) is a consequence of (\ref{e-1}) and the fact that
$$\lim_{n\to\infty}\frac{\sqrt{n}}{v_n}\sup_{t\in[0,1]}\left|\mathbb E\overline{\mathcal Q}_{3,t}^n(X^0)- \beta_{t}\right|=0,$$
which is an adaptation of the proof in Jiang \cite{Jiang1}.
\vspace{10pt}

\textbf{Proof  of Theorem \ref{MDP2}}
\vspace{10pt}

For (\ref{exponential-convergence-mdp}), we will prove that for $\ell=1,2$
$$
\frac{\sqrt{n}}{v_n}\sup_{t\in[0,1]}\left\|\mathcal{Q}_{\ell,t}^n(X^0)-Q_{\ell,t}^n\right\|\superexpmdp0 \quad{\rm and}\quad
\frac{\sqrt{n}}{v_n}\sup_{t\in[0,1]}\left\|\mathcal{C}_{t}^n(X^0)-C_{t}^n\right\|\superexpmdp0.
$$

From Lemma \ref{lem-mdp-1}, it follows that as $n\to\infty$
\begin{equation}\label{E-equv}
\frac{\sqrt{n}}{v_n}\sup_{t\in[0,1]}\left(\mathbb E(\mathcal{Q}_{\ell,t}^n(X^0)-Q_{\ell,t}^n)\vee
\mathbb E(\mathcal{C}_{t}^n(X^0)-C_{t}^n)\right)\to 0.
\end{equation}

Then, we only need to prove that
\begin{equation}\label{exponential-convergence-mdp-1}
\frac{\sqrt{n}}{v_n}\sup_{t\in[0,1]}\left\|\mathcal{Q}_{\ell,t}^n(X^0)-Q_{\ell,t}^n
-\mathbb E(\mathcal{Q}_{\ell,t}^n(X^0)-Q_{\ell,t}^n)\right\|\superexpmdp0
\end{equation}
and
\begin{equation}\label{exponential-convergence-mdp-2}
\frac{\sqrt{n}}{v_n}\sup_{t\in[0,1]}\left\|\mathcal{C}_{t}^n(X^0)-C_{t}^n-\mathbb E(\mathcal{C}_{t}^n(X^0)-C_{t}^n)\right\|\superexpmdp0.
\end{equation}

We start by the proof of (\ref{exponential-convergence-mdp-1}). Remark that
$\left(\mathcal{Q}_{\ell,t}^n(X^0)-Q_{\ell,t}^n-\mathbb E(\mathcal{Q}_{\ell,t}^n(X^0)-Q_{\ell,t}^n)\right)$
is a $\mathcal{F}_{[nt]/n}$-martingale. Then
$$
\exp\left(\lambda\left(\mathcal{Q}_{\ell,t}^n(X^0)-Q_{\ell,t}^n-\mathbb E(\mathcal{Q}_{\ell,t}^n(X^0)-Q_{\ell,t}^n)\right)\right)
$$
is a submartigale. By the maximal inequality, we have for any $\eta, \lambda>0$
\begin{align*}
&\mathbb{P}\left(\frac{\sqrt{n}}{v_n}\sup_{t\in[0,1]}\left(\mathcal{Q}_{\ell,t}^n(X^0)-Q_{\ell,t}^n
-\mathbb E(\mathcal{Q}_{\ell,t}^n(X^0)-Q_{\ell,t}^n)\right)>\eta\right)\\
&\leq e^{-\lambda v_n^2\eta}\mathbb{E}\exp\left(\lambda\sqrt{n}v_n
\left(\mathcal{Q}_{\ell,1}^n(X^0)-Q_{\ell,1}^n
-\mathbb E(\mathcal{Q}_{\ell,1}^n(X^0)-Q_{\ell,1}^n)\right)\right)
\end{align*}
and
\begin{align*}
&\mathbb{P}\left(\frac{\sqrt{n}}{v_n}\inf_{t\in[0,1]}\left(\mathcal{Q}_{\ell,t}^n(X^0)-Q_{\ell,t}^n
-\mathbb E(\mathcal{Q}_{\ell,t}^n(X^0)-Q_{\ell,t}^n)\right)<-\eta\right)\\
&\leq e^{-\lambda v_n^2\eta}\mathbb{E}\exp\left(-\lambda\sqrt{n}v_n
\left(\mathcal{Q}_{\ell,1}^n(X^0)-Q_{\ell,1}^n
-\mathbb E(\mathcal{Q}_{\ell,1}^n(X^0)-Q_{\ell,1}^n)\right)\right).
\end{align*}

Together with (\ref{E-e-1}) and (\ref{E-equv}), we have
$$
\limsup_{n\to\infty}\frac{1}{v_n^2}\log\mathbb{P}\left(\frac{\sqrt{n}}{v_n}
\sup_{t\in[0,1]}\left|\mathcal{Q}_{\ell,t}^n(X^0)-Q_{\ell,t}^n
-\mathbb E(\mathcal{Q}_{\ell,t}^n(X^0)-Q_{\ell,t}^n)\right|>\eta\right)\le -\lambda\eta.
$$
(\ref{exponential-convergence-mdp-1}) can be obtained by letting $\lambda$ goes to infinity.

Similarly, we can have (\ref{exponential-convergence-mdp-2}) by (\ref{E-e-1}), (\ref{decompC}) and (\ref{E-equv}).
\vspace{10pt}


\subsection{Proof of Theorem \ref{LDP1}}$\qquad$

We will do the proof in two steps.
\vspace{10pt}

{\underline {\it Step 1}}  We will prove that
$$
\mathcal V_1^n(X^0)-V_1^n\superexpldp0.
$$

For that, we will prove  that for $\ell =1,2$
\begin{equation}\label{nqldp}
\mathcal Q_{\ell,1}^n(X^0)-Q_{\ell,1}^n \superexpldp 0,
\end{equation}
and
\begin{equation}\label{ncldp}
\mathcal C_1^n(X^0)- C_1^n\superexpldp 0.
\end{equation}

We start by the proof of (\ref{nqldp}). Since the processes $X_{\ell}$ and $D_{\ell}$ have independent increment, by Chebyshev inequality we obtain for all $\theta>0$
$$
\mathbb P\left(\mathcal Q_{\ell,1}^n(X^0)-Q_{\ell,1}^n>\delta\right)
  \le e^{-\theta n\delta}\prod_{k=1}^{n}\mathbb E\left(e^{\theta n\left[(\Delta_k^nX^0_{\ell})^2{\bf 1}_{\{(\Delta_k^nX^0_{\ell})^2\le r(\frac 1n)\}}-(\Delta_k^nD_{\ell})^2\right]}\right).
$$

Similar to (\ref{Q}),
$$
\mathbb E\left(e^{\theta n\left[(\Delta_k^nX^0_{\ell})^2{\bf 1}_{\{(\Delta_k^nX^0_{\ell})^2\le r(\frac 1n)\}}
-(\Delta_k^nD_{\ell})^2\right]}\right) \le I_1(k,n)+I_2(k,n),$$
where
$$ I_1(k,n):=\mathbb E\left(e^{\theta n\left[(\Delta_k^nX^0_{\ell})^2-(\Delta_k^nD_{\ell})^2\right]}{\bf 1}_{\{(\Delta_k^nX^0_{\ell})^2\le r(\frac 1n)\}}\right)$$
and
$$I_2(k,n):=\mathbb P\left( (\Delta_k^nX^0_{\ell})^2> r(\frac 1n) \right)$$

From (\ref{r0}), (\ref{r1}) and (\ref{r2}), it follows that
$$
I_2(k,n)\le\exp\left(-\frac{r(\frac 1n)} {2\int_{t_{k-1}}^{t_k}\sigma_{\ell,s}^2ds}\right)+(1-e^{-\lambda_{\ell}/n}).
$$
and
\begin{eqnarray*}
I_1(k,n)\le 1+\mathbb E\left(e^{\theta n\left[(\Delta_k^nX^0_{\ell})^2-(\Delta_k^nD_{\ell})^2\right]}{\bf 1}_{\{(\Delta_k^nX^0_{\ell})^2\le r(\frac 1n), \Delta_k^nN_{\ell}\not=0 \}}\right).
\end{eqnarray*}

Let $(\alpha_n)$ be a sequence of real numbers such that  $\alpha_n\rightarrow 0$, which will be chosen latter. We  have
$$
\mathbb E\left(e^{\theta n(\Delta_k^nX^0_{\ell})^2}{\bf 1}_{\{(\Delta_k^nX^0_{\ell})^2\le r(\frac 1n), \Delta_k^nN_{\ell}\not=0 \}}\right)=F_1(k,n)+F_2(k,n),$$
where
$$F_1(k,n):=\mathbb E\left(e^{\theta n(\Delta_k^nX^0_{\ell})^2}{\bf 1}_{\{(\Delta_k^nX^0_{\ell})^2\le r(\frac 1n), \Delta_k^nN_{\ell}\not=0, | \Delta_k^nJ_{\ell}| \le \alpha_n\}}\right)$$
and
$$F_2(k,n):=\mathbb E\left(e^{\theta n(\Delta_k^nX^0_{\ell})^2}{\bf 1}_{\{(\Delta_k^nX^0_{\ell})^2\le r(\frac 1n), \Delta_k^nN_{\ell}\not=0, | \Delta_k^nJ_{\ell}| > \alpha_n\}}\right).$$

We have to prove that for $\ell=1,2$ $\lim_{n\rightarrow\infty}\max_{k=1}^nF_{\ell}(k,n)\rightarrow 0.$ We start with $F_2(k,n)$. 

From condition ({\bf LDP}), it follows that $n\max_{k=1}^{n} \int_{t_{k-1}}^{t_k} \sigma_{\ell,s}^2ds <+\infty$. 

So for all $\theta>0$, we choose
$$\alpha_n=\left(2\sqrt{\theta n\max_{k=1}^n\int_{t_{k-1}}^{t_k} \sigma_{\ell,s}^2ds}+1\right) \sqrt{r(1/n)}.$$

Then it is easy to see that
$$
F_2(k,n)\le e^{\theta n r(\frac 1n)}    \mathbb P  \left( |Z|\ge\frac{2\sqrt{\theta n\max_{k=1}^n\int_{t_{k-1}}^{t_k} \sigma_{\ell,s}^2ds}\sqrt{r(\frac 1n)}}{\sqrt {\int_{t_{k-1}}^{t_k} \sigma_{\ell,s}^2ds}} \right),
$$
where $Z$ is a standard Gaussian random variable. As a consequence of the well-known inequality $\int_y^{+\infty}e^{-\frac{z^2}{2}}dz\le (1/y)e^{-\frac{y^2}{2}}$, for all $y>0$, we obtain

 $$F_2(k,n)\le e^{\theta n r(\frac 1n)} \sqrt{\frac{2}{\pi}}\frac{1}{\sqrt{\theta nr(1/n)}}e^{-2\theta nr(\frac 1n)}.$$

So for $n$ large enough and $\theta>1$, we have
$$\max_{k=1}^nF_2(k,n)\le e^{-\theta n r(\frac 1n)} \longrightarrow 0\quad{\rm as}\quad n\rightarrow \infty.$$

Now we will control $F_1(k,n)$. Using the fact that
$$
\theta n (\Delta_k^nX^0_{\ell})^2\le\theta n\left[\frac{1}{4\theta n\max_{k=1}^n\int_{t_{k-1}}^{t_k}\sigma_{\ell,s}^2ds}(\Delta_k^nD_{\ell})^2+4\theta n\max_{k=1}^n\int_{t_{k-1}}^{t_k}\sigma_{\ell,s}^2ds(\Delta_k^nJ_{\ell})^2\right],
$$
we  have with the same choose of the sequence $\alpha_n$, by independence of  $\Delta_k^nD_{\ell}$ and $\Delta_k^nJ_{\ell}$ and Cauchy-Schwarz inequality that
\begin{eqnarray*}
F_1(k,n)&\le& \mathbb E\left(e^{\frac{(\Delta_k^nD_{\ell})^2}{4\max_{k=1}^n\int_{t_{k-1}}^{t_k}\sigma_{\ell,s}^2ds}} \right) \mathbb E\left(  
e^{4\theta^2 \left(n\max_{k=1}^n\int_{t_{k-1}}^{t_k}\sigma_{\ell,s}^2ds\right) n (\Delta_k^nJ_{\ell})^2}  {\bf 1}_{\{| \Delta_k^nJ_{\ell}| \le \alpha_n\}}
{\bf 1}_{\{ \Delta_k^nN_{\ell}\not=0\}}\right)\\
&\le&  \mathbb E\left(e^{\frac{Z^2}{4}} \right) \mathbb E^{\frac 12}\left(e^{8\theta^2 n(\Delta_k^nJ_{\ell})^2} {\bf 1}_{\{| \Delta_k^nJ_{\ell}| \le \alpha_n\}} \right)
\mathbb P^{\frac 12}\left( \Delta_k^nN_{\ell}\not=0\right).
\end{eqnarray*}

From Mancini \cite{Mancini1} page 877, we conclude that
$$\lim_{n\rightarrow \infty}\max_{k=1}^n\mathbb E\left(e^{8\theta^2 n(\Delta_k^nJ_{\ell})^2} {\bf 1}_{\{| \Delta_k^nJ_{\ell}| \le \alpha_n\}} \right)<\infty.$$

Since $Z$ is a standard Gaussian random variable, we conclude that
$$E\left(e^{\frac{Z^2}{4}} \right) <\infty.$$

 So that  $\max_{k=1}^nF_1(k,n)\le C (1-e^{-\lambda_{\ell}/n})\longrightarrow 0$ as $n\rightarrow \infty.$

Therefore,
$$
\lim_{n\to\infty}\frac{1}{n}\log\prod_{k=1}^{n}\mathbb E\left(e^{\theta n\left[(\Delta_k^nX^0_{\ell})^2{\bf 1}_{\{(\Delta_k^nX^0_{\ell})^2\le r(\frac 1n)\}}-(\Delta_k^nD_{\ell})^2\right]}\right)=0,
$$
which implies that for any $\theta>1$
$$
\lim_{n\to\infty}\frac{1}{n}\log\mathbb P\left(\mathcal Q_{\ell,1}^n(X^0)-Q_{\ell,1}^n>\delta\right)\leq-\theta\delta.
$$
Letting $\theta$ goes to infinity, we obtain that the left term in the last inequality goes to $-\infty$.
And similarly, by doing the same calculation with
$$
\mathbb P\left(\mathcal Q_{\ell,1}^n(X^0)-Q_{\ell,1}^n<-\delta\right),
$$
we can get (\ref{nqldp}).

To prove (\ref{ncldp}), we use the decomposition (\ref{decompC}) and an adaptation of the proof of  (\ref{nqldp}).

\vspace{10pt}

{\underline {\it Step 2}}  We will prove that
$$
\mathcal V_1^n(X)-\mathcal V_1^n(X^0)\superexpldp0.
$$

For that we use (\ref{QQ}) and (\ref{C}) and we choose $\varepsilon(n)$ such that $n\varepsilon(n)\rightarrow 0$ to obtain the result.
\vspace{10pt}


\subsection{Proof of Theorem \ref{LDP2}}$\qquad$

We will prove that for $\ell=1,2$
$$
\sup_{t\in[0,1]}\left\|\mathcal{Q}_{\ell,t}^n(X^0)-Q_{\ell,t}^n\right\|\superexpldp0\quad{\rm and}\quad
\sup_{t\in[0,1]}\left\|\mathcal{C}_{t}^n(X^0)-C_{t}^n\right\|\superexpldp0.
$$

To do that we use the same argument as in the proof of Theorem \ref{MDP2} and the fact that
$$\sup_{t\in[0,1]}\left |\mathbb E(\mathcal{Q}_{\ell,t}^n(X^0)-Q_{\ell,t}^n))\right| \longrightarrow 0.$$


\medskip

\nocite{*}

\bibliographystyle{acm}
\bibliography{bibliosoumis}

\nocite{*}

\vspace{10pt}

\end{document}